\documentclass[a4paper,10pt,reqno]{amsart}
\usepackage[a4paper,tmargin=40mm,bmargin=35mm,hmargin=30mm]{geometry}
\usepackage[T1]{fontenc}
\usepackage{lmodern}
\usepackage{amsmath}
\usepackage{amssymb}
\usepackage{amsthm}
\usepackage{stmaryrd}
\usepackage{tikz}
\usepackage[all]{xy}
\usepackage{booktabs}
\usepackage{multirow}
\usepackage{multicol}
\usepackage{here}
\usepackage{aliascnt}
\usepackage{hyperref}
\usepackage[noabbrev]{cleveref}

\newcommand{\NewTheorem}[2]{
	\newaliascnt{#1}{TheoremEnvironment}
	\newtheorem{#1}[#1]{#1}
	\aliascntresetthe{#1}
	\crefname{#1}{#1}{#2}
	\Crefname{#1}{#1}{#2}
}

\theoremstyle{definition}

\NewTheorem{Definition}{Definitions}
\NewTheorem{Remark}{Remarks}
\NewTheorem{Axiom}{Axioms}
\NewTheorem{Example}{Examples}
\NewTheorem{Observation}{Observations}
\NewTheorem{Convention}{Conventions}
\NewTheorem{Notation}{Notations}
\NewTheorem{Setting}{Settings}
\NewTheorem{Question}{Questions}
\NewTheorem{Answer}{Answers}
\NewTheorem{Conjecture}{Conjectures}
\NewTheorem{Problem}{Problems}
\NewTheorem{Solution}{Solutions}
\NewTheorem{Goal}{Goals}
\NewTheorem{Comment}{Comments}
\NewTheorem{Aim}{Aims}
\NewTheorem{Caution}{Cautions}
\NewTheorem{Exercise}{Exercises}
\NewTheorem{Examples}{Examples}
\theoremstyle{plain}
\NewTheorem{Proposition}{Propositions}
\NewTheorem{Lemma}{Lemmas}
\NewTheorem{Theorem}{Theorems}
\NewTheorem{Corollary}{Corollaries}

\crefname{enumi}{}{}
\Crefname{enumi}{}{}
\creflabelformat{enumi}{(#2#1#3)}
\crefname{enumii}{}{}
\Crefname{enumii}{}{}
\creflabelformat{enumii}{(#2#1#3)}
\crefname{enumiii}{}{}
\Crefname{enumiii}{}{}
\creflabelformat{enumiii}{(#2#1#3)}

\makeatletter
\renewcommand{\p@enumii}{}
\renewcommand{\p@enumiii}{}
\makeatother

\numberwithin{equation}{section}
\crefname{equation}{}{}
\Crefname{equation}{}{}
\creflabelformat{equation}{(#2#1#3)}

\newcommand{\SwapSymbols}[1]{
	\expandafter\let\expandafter\temporarysymbol\csname #1\endcsname
	\expandafter\let\csname #1\expandafter\endcsname\csname var#1\endcsname
	\expandafter\let\csname var#1\endcsname\temporarysymbol
}

\SwapSymbols{epsilon}
\SwapSymbols{phi}
\SwapSymbols{Gamma}
\SwapSymbols{Delta}
\SwapSymbols{Theta}
\SwapSymbols{Lambda}
\SwapSymbols{Xi}
\SwapSymbols{Pi}
\SwapSymbols{Sigma}
\SwapSymbols{Upsilon}
\SwapSymbols{Phi}
\SwapSymbols{Psi}
\SwapSymbols{Omega}

\newcommand{\cS}{\mathcal{S}}

\newcommand{\Modr}{\mathrm{Mod}\text{-}}

\newcommand{\To}{\longrightarrow}

\DeclareMathOperator{\Hom}{Hom}

\DeclareMathOperator{\Ext}{Ext}

\DeclareMathOperator{\Ann}{Ann}

\DeclareMathOperator{\Ker}{Ker}
\DeclareMathOperator{\Coker}{Coker}
\let\Im\relax
\DeclareMathOperator{\Im}{Im}

\DeclareMathOperator{\Supp}{Supp}

\DeclareMathOperator{\depth}{depth}

\title{Cofiniteness of modules and local cohomology}

\author{Hajar Sabzeh and Reza Sazeedeh}
\address{Department of Mathematics, Urmia University, P.O.Box: 165, Urmia, Iran}
\email{h.sabzeh96@gmail.com}
\address{Department of Mathematics, Urmia University, P.O.Box: 165, Urmia, Iran}
\email{rsazeedeh@ipm.ir}

\begin{document}

\begin{abstract}
Let $A$ be a commutative noetherian ring, let $\mathfrak a$ be an ideal of $A$ and let $n$ be a non-negative integer. In this paper, we study $\cS_n(\frak a)$, a certain class of $A$-modules and we find some sufficient conditions so that a module belongs to $\cS_n(\frak a)$. Moreover, we  study the cofiniteness of local cohomology modules when $\dim A/\frak a\geq 3$.
\end{abstract}

\maketitle
\tableofcontents

\section{Introduction}
Throughout this paper, assume that $A$ is a commutative noetherian ring, $\frak a$ is an ideal of $A$, $M$ is an $A$-module and $n$ is a non-negative integer. Following [H], the $A$-module $M$ is said to be $\frak a$-{\it cofinite} if $\Supp_A(M)\subseteq V(\frak a)$ and
$\Ext_A^i(A/\frak a, M)$ are finitely generated for all integers $i\geq 0$. The authors [KS] studied a criterion for cofinitenss of modules. We denote by $\cS_n(\frak a)$ the class of all $A$-modules $M$ satisfying the  following implication:
\begin{center}
   If $\Ext_A^i(A/\frak a, M)$ is
finitely generated for all $i\leq n$ and $\Supp(M)\subseteq V(\frak a)$,
then $M$ is $\frak a$-cofinite.
\end{center}

The abelianess of the category of $\frak a$-cofinite modules is of interest to many mathematicians  working in commutative algebra. This subject has been studied for small dimensions by various authors (cf. [K, M1, M2, BNS1]). As the Koszul complexes are  effective tools for going down the dimensions, Sazeedeh [S] studied the cofiniteness of Koszul cohomology modules.

 In Section 2, by means of the Koszul cohomologies of a module, we want to find out when this module belongs to $\cS_n(\frak a)$. Let ${\bf x}=x_1,\dots,x_t$ be a sequence of elements of $\frak a$. We show that if $\Ext_A^i(A/\frak a,M)$ is finitely generated for all $i\leq n+1$ and $H^i({\bf x},M)\in\cS_1(\frak a)$ for all $i\leq n$, then $H^i({\bf x},M)$ is $\frak a$-cofinite for all $i\leq n $ (cf. \cref{s1}). Moreover, we prove the following theorem.
 \begin{Theorem}
 If $H^i({\bf x},M)\in\cS_1(\frak a)$ for all $i\geq 0$, then $M\in\cS_{t+1}(\frak a)$.
\end{Theorem}
As an application of this theorem, if $\dim H^i({\bf x},M)\leq 1$ for all $i\geq 0$, then $M\in\cS_{t+1}(\frak a)$ (cf. \cref{c1}). Let $\Ext_A^i(A/\frak a,M)$ is finitely generated for all $i\leq n+2$ and let $H^i({\bf x},M)\in\cS_2(\frak a)$ for all $i\leq n$. Then we show that
 $H^i({\bf x},M)$ is $\frak a$-cofinite for all $i\leq n $ if and only if $\Hom_A(A/\frak a,H^i({\bf x},M))$ is finitely generated for all $i\leq n+1$ (cf. \cref{s2}). Moreover, we have the following theorem.
 \begin{Theorem}
Let $H^i({\bf x},M)\in\cS_2(\frak a)$ such that $\Hom_A(A/\frak a,H^i({\bf x},M))$ is finitely generated for all $i\geq 0$. Then $M\in\cS_{t+2}(\frak a)$.
\end{Theorem}
   Let $\Ext_A^i(A/\frak a,M)$ be finitely generated for all $i\leq n+3$ and let $H^i({\bf x},M)\in\cS_3(\frak a)$ for all $i\leq n$. Then we prove that if $H^i({\bf x},M)$ is $\frak a$-cofinite for all $i\leq n$, then $\Ext_A^j(A/\frak a,H^i({\bf x},M))$ is finitely generated for all $i\leq n+1$; $j=0,1$; furthermore if $\Hom_A(A/\frak a,H^{n+2}({\bf x},M))$ is finitely generated, then the converse also holds (cf. \cref{s3}). Moreover, we have the following theorem.
    \begin{Theorem}
Let $H^i({\bf x},M)\in\cS_3(\frak a)$ such that $\Hom_A(A/\frak a,H^i({\bf x},M))$ and $\Ext_A^1(A/\frak a,H^i({\bf x},M))$ are finitely generated for all $i\geq 0$. Then $M\in\cS_{t+3}(\frak a)$.
\end{Theorem}
  
 As an application of this theorem, let $A$ be a local ring such that $\dim A/({\bf x})+xA)\leq 3$ for some $x\in\frak a$. If
$\Hom_A(A/\frak a,H^i({\bf x},M)/xH^i({\bf x},M))$ and $\Ext_A^j(A/\frak a,H^i({\bf x},M))$ for all  $i\geq 0$ and $j=0,1$ are finitely generated, then $M\in\cS_{t+3}(\frak a)$ (cf. \cref{co3}).

   In Section 3, we study the cofiniteness of local cohomology modules. Let $\Ext^i_A(A/\frak a,M)$ be finitely generated for all $i\geq 0$ and let $t<s$ be non-negative integers such that $H_{\frak a}^i(M)$ is $\frak a$-cofinite for all $i\neq t,s$. Then we show that  $H_{\frak a}^{t}(M)\in\cS_{n+s-t+1}(\frak a)$ if and only if $H_{\frak a}^s(M)\in\cS_{n}(\frak a)$ (cf. \cref{tt}). Moreover, we have the following theorem.
 \begin{Theorem}
 Let $\dim A/\frak a=d\geq 3$ and let $\depth(\Ann(M),A/\frak a)\geq d-2$. If $\Ext_A^i(A/\frak a,M)$ is finitely generated for all $i\leq n+1$, then $\Hom_A(A/\frak a,H_{\frak a}^i(M))$ is finitely generated for all $i\leq n$ if and only if $H_{\frak a}^i(M)$ is $\frak a$-cofinite for all $i<n.$ 
 \end{Theorem}
  In the end of this paper we get a similar result for those rings $A$ which $\dim A\geq 4$.

\section{A criterion for cofiniteness of  modules}
 Throughout this section, $M$ is an $A$-module with $\Supp_AM\subseteq V(\frak a)$ and $n$ is a non-negative integer and ${\bf x}=x_1,\dots, x_t$ is a sequence of elements of $\frak a$.
 
 An $A$-module $M$ is said to be $\frak a$-{\it cofinite} if $\Supp_A(M)\subseteq V(\frak a)$ and
$\Ext_A^i(A/\frak a, M)$ are finitely generated for all integers $i\geq 0$.

For a non-negative integer $n$, we denote by $\cS_n(\frak a)$ the class of all $A$-modules $M$ satisfying the  following implication:
\begin{center}
   If $\Ext_A^i(A/\frak a, M)$ is
finitely generated for all $i\leq n$ and $\Supp(M)\subseteq V(\frak a)$,
then $M$ is $\frak a$-cofinite.
\end{center}

\begin{Examples}\label{ex}
(i) Assume that $\frak a$ is an arbitrary ideal of $A$ and $M$ is an $A$-module of dimension $d$ where $\dim M$ means the dimension of $\Supp_A M$; which is the length of the longest chain of prime ideals in $\Supp_AM$. Then $H_{\frak a}^d(M)$ is in $\cS_0(\frak a)$. To be more precise, if $\Hom_A(A/\frak a,H_{\frak a}^d(M))$ is a finitely generated $A$-module, then it follows from \cite[Theorem 3.11]{NS} that $H_{\frak a}^d(M)$ is artinian and so, since $\Hom_A(A/\frak a,H_{\frak a}^d(M))$ has finite length, according to \cite[Proposition 4.1]{M2}, the module $H_{\frak a}^d(M)$ is $\frak a$-cofinite.
  
(ii) Given an arbitrary ideal $\frak a$ of $A$, by virtue of \cite[Proposition 2.6]{BNS1}, $M\in\cS_1(\frak a)$ for all modules $M$ with $\dim M\leq 1$. Especially, if $\dim A/\frak a=1$, then it follows from \cite[Theorem 2.3]{M2} that $\cS_1(\frak a)=\Modr A$, the category of all $A$-modules. Furthermore, if $\dim A=2$, then it follows from \cite[Corollary 2.4]{NS} that $\cS_1(\frak a)=\Modr A$ for any ideal $\frak a$ of $A$. 

(iii) Let $\frak a$ be an ideal of a local ring $A$ with $\dim A/\frak a=2$. It follows from \cite[Theorem 3.5]{BNS2} that $\cS_2(\frak a)=\Modr A$. Furthermore, if $A$ is a local ring with $\dim A=3$, then it follows from \cite[Corollary 2.5]{NS} that $\cS_2(\frak a)=\Modr A$ for any ideal $\frak a$ of $A$. 
\end{Examples}

\medskip
When the Koszul cohomology modules of an $A$-module belong to $\cS_1(\frak a)$, the following result shows that they are $\frak a$-cofinite.
\begin{Proposition}\label{s1}
Let $\Ext_A^i(A/\frak a,M)$ is finitely generated for all $i\leq n+1$. If $H^i({\bf x},M)\in\cS_1(\frak a)$ for all $i\leq n$, then $H^i({\bf x},M)$ is $\frak a$-cofinite for all $i\leq n $.
\end{Proposition}
\begin{proof}
Consider the Koszul complex $$K^*({\bf x},M):0\To K^0\stackrel{d^0}\To K^1\stackrel{d^1}\To\dots\stackrel{d^{t-1}}\To K^t\To 0 $$
and assume that $Z^i=\Ker d^i$, $B^i=\Im d^{i-1}$, $C^i=\Coker d^i$ and $H^i=H^i({\bf x},M)$ for each $i$. Then for each $j\geq 0$, we have an exact sequence of modules

$$0\To H^j\To C^j\To K^{j+1}\To C^{j+1}\To o\hspace{1cm}(\dag_j).$$
We prove by induction on $i$ that $H^i$ is $\frak a$-cofinite and $\Ext_A^j(A/\frak a,C^{i+1})$ is finitely generated for all $i\leq n$ and all $j\leq n-i$.
 Assume that $i<n$ and so the induction hypothesis implies that $\Ext_A^j(A/\frak a,C^i)$ is finitely generated for all $j\leq n+1-i$. In view of $(\dag_i)$ it is clear that $\Ext_A^j(A/\frak a,H^i)$ is finitely generated for $i=0,1$ and so the fact that $H^i\in\cS_1(\frak a)$ forces $H^i$ is $\frak a$-cofinite and $\Ext_A^j(A/\frak a,C^{i+1})$ is finitely generated  for all $j\leq n-i$. 
\end{proof}
\medskip

The following theorem provides a sufficient condition so that a module belongs to $\cS_{t+1}(\frak a)$.

\begin{Theorem}\label{ts1}
If $H^i({\bf x},M)\in\cS_1(\frak a)$ for all $i\geq 0$, then $M\in\cS_{t+1}(\frak a)$.
\end{Theorem}
\begin{proof}
Assume that $\Ext_A^i(A/\frak a,M)$ is finitely generated for all $i\leq t+1$. Then, it follows from \cref{s1} that $H^i({\bf x},M)$ is $\frak a$-cofinite for all $i\geq 0$. Consequently, [S, Theorem 2.4] implies that $M$ is $\frak a$-cofinite.   
\end{proof}

\medskip

\begin{Corollary}\label{c1}
If $\dim H^i({\bf x},M)\leq 1$ for all $i\geq 0$, then $M\in\cS_{t+1}(\frak a)$.
\end{Corollary}
\begin{proof}
By virtue of \cref{ex}, $H^i({\bf x},M)\in\cS_1(\frak a)$ for all $i\geq 0$ and so the result follows from \cref{ts1}.  
\end{proof}

\medskip
When the Koszul cohomology modules of an $A$-module belong to $\cS_2(\frak a)$, we have the following result about their $\frak a$-cofiniteness .

\begin{Proposition}\label{s2}
Let $\Ext_A^i(A/\frak a, M)$ is finitely generated for all $i\leq n+2$ and let $H^i({\bf x},M)\in\cS_2(\frak a)$ for all $i\leq n$. Then the following conditions hold.

${\rm (i)}$ $H^i({\bf x},M)$ is $\frak a$-cofinite for all $i\leq n $.

${\rm (ii)}$ $\Hom_A(A/\frak a,H^i({\bf x},M))$ is finitely generated for all $i\leq n+1$. 
\end{Proposition}
\begin{proof}
Consider the same notation as in the proof of \cref{s1}. (i)$\Rightarrow$(ii). By the assumption $\Hom_A(A/\frak a, H^i)$ is finitely generated for all $i\leq n$. Thus, it is straightforward to see that $\Ext_A^j(A/\frak a, C^i)$ is finitely generated for all $i\leq n+1$ and all $j\leq n+2-i$.  Since $\Hom_A(A/\frak a,C^{n+1})$ is finitely generated, $(\dag_{n+1})$ implies that $\Hom_A(A/\frak a,H^{n+1})$ is finitely generated.   
(ii)$\Rightarrow$(i). We prove by induction on $i$ that $H^i$ is $\frak a$-cofinite and $\Ext_A^j((A/\frak a,C^{i+1})$ is finitely generated for all $i\leq n$ and all $j\leq n+1-i$. The induction hypothesis implies that  $\Ext_A^j((A/\frak a,C^{i})$ is finitely generated for all $j\leq n+2-i$. Since by the assumption $\Hom_A(A/\frak a,H^{i+1})$ is finitely generated, $(\dag_{i+1})$ implies that $\Hom_A(A/\frak a,C^{i+1})$ is finitely generated; and hence $(\dag_i)$ implies that $\Ext_A^j(A/\frak a,H^i)$ is finitely generated for all $j\leq 2$. Now, since $H^i\in\cS_2(\frak a)$, it is $\frak a$-cofinite; furthermore $(\dag_i)$ implies that $\Ext_A^j((A/\frak a,C^{i+1})$ is finitely generated for all $i\leq n$ and all $j\leq n+1-i$. 
\end{proof}

\medskip

The following theorem provides a  sufficient condition so that a module belongs to $\cS_{t+2}(\frak a)$.
\begin{Theorem}\label{ts2}
Let $H^i({\bf x},M)\in\cS_2(\frak a)$ such that $\Hom_A(A/\frak a,H^i({\bf x},M))$ is finitely generated for all $i\geq 0$. Then $M\in\cS_{t+2}(\frak a)$.
\end{Theorem}
\begin{proof}
Assume that $\Ext_A^i(A/\frak a,M)$ is finitely generated for all $i\leq t+2$. Then it follows from \cref{s2} that $H^i({\bf x},M)$ is $\frak a$-cofinite for all $i\geq 0$. Consequently, [S, Theorem 2.4] implies that $M$ is $\frak a$-cofinite. 
\end{proof}

\medskip

\begin{Corollary}
Let $A$ be a local ring such that $\dim A/({\bf x})\leq 3$ and let $\Hom_A(A/\frak a,H^i({\bf x},M))$ be finitely generated for all $i\geq 0$. Then $M\in\cS_{t+2}(\frak a)$.  
\end{Corollary}
\begin{proof}
 Put $B=A/({\bf x})$ and $\frak b=\frak aB$. We observe that $H^i({\bf x},M)$ is a $B$-module for each $i\geq 0$ and since $\dim B\leq 3$, it follows from \cref{ex} (iii) that $H^i({\bf x},M)\in\cS_2(\frak aB)$ for each $i\geq 0$. Now [KS, Proposition 2.15] implies that $H^i({\bf x},M)\in\cS_2(\frak a)$; and consequently, it follows from \cref{ts2} that $M\in\cS_{t+2}(\frak a)$.
\end{proof}

\medskip

When the Koszul cohomology modules of an $A$-module belong to $\cS_3(\frak a)$, we have the following result about their $\frak a$-cofiniteness

\begin{Proposition}\label{s3}
 Let $\Ext_A^i(A/\frak a, M)$ be finitely generated for all $i\leq n+3$ and let $H^i({\bf x},M)\in\cS_3(\frak a)$ for all $i\leq n$. Consider the following statements.

${\rm (i)}$ $H^i({\bf x},M)$ is $\frak a$-cofinite for all $i\leq n $.

${\rm (ii)}$ $\Ext_A^j(A/\frak a,H^i({\bf x},M))$ is finitely generated for $j=0, 1$ and all $i\leq n+1.$

 Then ${\rm (i)}\Rightarrow{\rm (ii)}$ holds. Moreover, if $\Hom_A(A/\frak a,H^{n+2}({\bf x},M))$ is finitely generated, then ${\rm (ii)}\Rightarrow{\rm (i)}$ holds. 
\end{Proposition}
\begin{proof}
Consider the same notation as in the proof of \cref{s1}. (i)$\Rightarrow$(ii). Clearly $\Ext_A^j(A/\frak a,H^i)$ is finitely generated for all $i\leq n$ and $j=0,1$; furthermore $\Ext_A^j(A/\frak a, C^i)$ is finitely generated for all $i\leq n+1$ and all $j\leq n+3-i$.  Since $\Ext_A^j(A/\frak a,C^{n+1})$ is finitely generated for $j=0,1$, the exact sequence $(\dag_{n+1})$ implies that $\Ext_A^j(A/\frak a,H^{n+1})$ is finitely generated for $j=0,1$. (ii)$\Rightarrow$(i). We prove by induction on $i$ that $H^i$ is $\frak a$-cofinite and $\Ext_A^j(A/\frak a,C^{i+1})$ is finitely generated for all $i\leq n$ and all $j\leq n+2-i$. The exact sequence $(\dag_{i+2})$ implies that $\Hom_A(A/\frak a,C^{i+2})$ is finitely generated and so it follows from $(\dag_{i+1})$ that $\Ext^j(A/\frak a, C^{i+1})$ is finitely generated for $j=0,1$. Since by the induction hypothesis, $\Ext_A^j(A/\frak a,C^i)$ is finitely generated  for all $j\leq n+3-i$, the exact sequence $(\dag_i)$ implies that $\Ext_A^j(A/\frak a,H^i)$ is finitely generated for all $j\leq n+3-i$; and hence since $H^i\in\cS_3(\frak a)$, we deduce that $H^i$ is $\frak a$-cofinite.  
\end{proof}

\medskip

The following theorem provides a  sufficient condition so that a module belongs to $\cS_{t+3}(\frak a)$.

\begin{Theorem}\label{ts3}
Let $H^i({\bf x},M)\in\cS_3(\frak a)$ such that $\Hom_A(A/\frak a,H^i({\bf x},M))$ and $\Ext_A^1(A/\frak a,H^i({\bf x},M))$ are finitely generated for all $i\geq 0$. Then $M\in\cS_{t+3}(\frak a)$.
\end{Theorem}
\begin{proof}
Assume that $\Ext_A^i(A/\frak a,M)$ is finitely generated for all $i\leq t+3$. Then it follows from \cref{s3} that $H^i({\bf x},M)$ is $\frak a$-cofinite for all $i\geq 0$. Consequently, [S, Theorem 2.4] implies that $M$ is $\frak a$-cofinite. 
\end{proof}

\medskip
\begin{Corollary}\label{co3}
Let $A$ be a local ring such that $\dim A/({\bf x})+xA)\leq 3$ for some $x\in\frak a$. If
$\Hom_A(A/\frak a,H^i({\bf x},M)/xH^i({\bf x},M))$ and $\Ext_A^j(A/\frak a,H^i({\bf x},M))$ for all  $i\geq 0$ and $j=0,1$ are finitely generated, then $M\in\cS_{t+3}(\frak a)$.  
\end{Corollary}
\begin{proof}
Set $B=A/xA$, $\frak b=\frak a/xA$. In view of the exact sequence $$H^i({\bf x},x,M)\To H^i({\bf x},M)\stackrel{x.}\To H^i({\bf x},M)\To H^{i+1}({\bf x},x,M)$$ for $i\geq 0$, there is an exact sequence $$0\To (0:_{H^i({\bf x},M)}x)\To H^i({\bf x},M)\stackrel{x.}\To H^i({\bf x},M)\To H^i({\bf x},M)/xH^i({\bf x},M)\To 0.$$
As $(0:_{H^i({\bf x},M)}x)$ and $H^i({\bf x},M)/xH^i({\bf x},M)$ are $B$-modules, they belong to $\cS_2(\frak b)$ by \cref{ex} (iii) and so 
by virtue of [KS, Proposition 2.15], they belong to $\cS_2(\frak a)$. Thus, \cref{ts2} implies that $H^i({\bf x},M)\in\cS_3(\frak a)$ for each $i\geq 0$ and consequently, $M\in\cS_{t+3}(\frak a)$ by using \cref{ts3}. 
\end{proof}


\section{Cofiniteness of local cohomology modules}
 Throughout this section, $M$ is an $A$-module, $\frak a$ is an ideal of $A$ and $n$ is a positive integer.
 We study the cofiniteness of local cohomology.For more details about local cohomology, we refer the reader to the textbook of Brodmann and Sharp [BS]. Throughout this section $n$ is a non-negative integer. 

\medskip
\begin{Theorem}\label{tt}
Let $\Ext^i_A(A/\frak a,M)$ be finitely generated for all $i\geq 0$ and let $t<s$ be non-negative integers such that $H_{\frak a}^i(M)$ is $\frak a$-cofinite for all $i\neq t,s$. Then $H_{\frak a}^{t}(M)\in\cS_{n+s-t+1}(\frak a)$ if and only if $H_{\frak a}^s(M)\in\cS_{n}(\frak a)$. 
\end{Theorem}
\begin{proof}
We first assume $H_{\frak a}^{t}(M)\in\cS_{n+s-t+1}(\frak a)$ and that  $\Ext_A^p(A/\frak a,H_{\frak a}^s(M))\in\cS_{n}(\frak a)$ is finitely generated for all $p\leq n$. There is the Grothendieck spectral sequence $$E_2^{p,q}:=\Ext_A^p(A/\frak a, H_{\frak a}^q(M))\Rightarrow \Ext_A^{p+q}(A/\frak a,M).$$
For each $r\geq 3$, consider the sequence $E_{r-1}^{p-r+1,t+r-2}\stackrel{d_{r-1}^{p-r+1,t+r+2}}\longrightarrow E_{r-1}^{p,t}\stackrel{d_{r-1}^{p,t}}\To E_{r-1}^{p+r-1,t-r+2}$ and so $E_r^{p,t}=\Ker d_{r-1}^{p,t}/\Im d_{r-1}^{p-r+1,t+r-2}$. Considering $p\leq n+s-t+1$, we show that $E_2^{p,t}$ is finitely generated . If $r=s-t+2$, then $p-r+1\leq n$ and so $E_{r-1}^{p-r+1,t+r-2}$ is finitely generated by the argument in the beginning of proof. If $r\neq s-t+2$, then the assumption implies that $E_{r-1}^{p-r+1,t+r-2}$ is finitely generated for all $p\geq 0$ (we observe that $t+r-2\neq t$). Consequently,  $\Im d_{r-1}^{p-r+1,t+r-2}$ is finitely generated for all $r\geq 3$ and $p\leq n+s-t+1$. But  there is a finite filtration 
$$0=\Phi^{p+t+1}H^{p+t}\subset \Phi^{p+t}H^{p+t}\subset\dots\subset\Phi^1H^{p+t}\subset \Phi^0H^{p+t}\subset H^{p+t}$$ 
such that $\Phi^pH^{p+t}/\Phi^{p+1}H^{p+t}\cong E_{\infty}^{p,t}$ for all $p\geq 0$. In view of the assumption, $H^{p+t}=\Ext_A^{p+t}(A/\frak a,M)$ is finitely generated and so $E_{\infty}^{p,t}$ is finitely generated for all $p\geq 0$. On the other hand, $E_r^{p,t}=E_{\infty}^{p,t}$ for sufficiently large $r$ and so $E_r^{p,t}$ is finitely generated for all $p\geq 0$. The previous argument implies that $\Ker d_{r-1}^{p,t}$ is finitely generated for all $p\leq n+s-t+1$; moreover since $r\geq 3$, we have $t-r+2\leq t-1$ and so the assumption implies that $E_{r-1}^{p+r-1,t-r+2}$ is finitely generated, and hence $E_{r-1}^{p,t}$ is finitely generated for all $p\leq n+s-t+1$. Continuing this way, we deduce that  $E_2^{p,t}$ is finitely generated for all $p\leq n+s-t+1$, and since $H_{\frak a}^t(M)\in\cS_{n+s-t+1}(\frak a)$, we deduce that $H_{\frak a}^t(M)$ is $\frak a$-cofinite. Therefore, $E_2^{p,q}$ is finitely generated for all $q\neq s$ and all $p\geq 0$. By a similar argument, we have $E_{\infty}^{p,s}=E_r^{p,s}$ for sufficiently large $r$ and all $p\geq 0$ and since $E_{\infty}^{p,s}$ is a subquotient of $\Ext_A^{p+s}(A/\frak a,M)$, it is finitely generated so that $E_r^{p,s}$ is finitely generated. In view of the sequence $E_{r-1}^{p-r+1,s+r-2}\stackrel{d_{r-1}^{p-r+1,s+r+2}}\longrightarrow E_{r-1}^{p,s}\stackrel{d_{r-1}^{p,s}}\To E_{r-1}^{p+r-1,s-r+2}$, since $E_{r-1}^{p-r+1,s+r-2}$ is finitely generated and $E_r^{p,s}=\Ker d_{r-1}^{p,s}/\Im d_{r-1}^{p-r+1,s+r-2}$, we conclude that  
$\Ker d_{r-1}^{p,s}$ is finitely generated; and hence $E_{r-1}^{p,s}$ is finitely generated. Continuing this way, we deduce that $E_3^{p,s}$ is finitely generated and so is $\Ker d_2^{p,s}$ for all $p\geq 0$. Now the exact sequence $0\To \Ker d_2^{p,s}\To E_2^{p,s}\To E_2^{p+2,s+1}$ implies that $E_2^{p,s}$ is finitely generated for all $p\geq 0$; and consequently $H_{\frak a}^s(M)$ is $\frak a$-cofinite. A similar proof gets the converse.
\end{proof}

The following lemma extends [KS, Proposition 2.6].

\medskip
\begin{Lemma}\label{three}
Let $\dim A/\frak a=3$ and $\depth(\Ann(M),A/\frak a)>0$. If $\Ext_A^i(A/\frak a,M)$ is finitely generated for all $i\leq n+1$, then the following conditions are equivalent.

${\rm (i)}$ $\Hom_A(A/\frak a,H_{\frak a}^i(M))$ is finitely generated for all $i\leq n$.

${\rm (ii)}$ $H_{\frak a}^i(M)$ is $\frak a$-cofinite for all $i<n$.  
\end{Lemma}
\begin{proof}
Since $\depth(\Ann(M),A/\frak a)>0$, there exists an element $x\in\Ann(M)$ such that $x$ is a non-zerodivisor on  $A/\frak a$. Taking $\frak b=\frak a+xA$, we have $\dim A/b=2$ and it follows from [DM, Proposition 1] that $\Ext_A^i(A/\frak b,M)$ is finitely generated for all $i\leq n+1$. Moreover, we have $\Gamma_{xA}(M)=\Gamma_{x\frak a}(M)=M$ and $\Gamma_{\frak a}(M)=\Gamma_{\frak b}(M)$. Thus the ideals $\frak a$ and $xA$ of $A$  provides the following Mayer-Vietoris exact sequence 
$$0\To \Gamma_{\frak b}(M)\To \Gamma_{\frak a}(M)\oplus M\To M\To H_{\frak b}^1(M)\To H_{\frak a}^1(M)\To 0$$ and the isomorphism  $H_{\frak a}^i(M)\cong H_{\frak b}^i(M)$ for each $i\geq 2$. (i)$\Rightarrow$(ii). The case $n=1$ follows from [KS, Proposition 2.6]. For $n\geq 2$,  $\Gamma_{\frak a}(M)$ is $\frak a$-cofinite. Then $\Gamma_{\frak b}(M)=\Gamma_{\frak a}(M)$ is $\frak b$-cofinite too. Moreover, it is clear that $\Hom_A(A/\frak b,H_{\frak a}^i(M))$ is finitely generated for all $i\leq n$. Applying the functor $\Hom_A(A/\frak b,-)$ to the above exact sequence, we deduce that $\Hom_A(A/\frak b,\Gamma_{\frak b}(M))$ and $\Hom_A(A/\frak b,H_{\frak b}^1(M))$ are finitely generated. Furthermore $\Hom_A(A/\frak b,H_{\frak b}^i(M))\cong\Hom_A(A/\frak b,H_{\frak a}^i(M))$ is finitely generated for each $2\leq i\leq n$. Now, [NS, Theorem 3.7] implies that $H_{\frak b}^i(M)$ is $\frak b$-cofinite and consequently $H_{\frak a}^i(M)$ is $\frak a$-cofinite for all $i<n$ using [DM, Proposition 2]. (ii)$\Rightarrow$(i) Since $H_{\frak a}^i(M)$ is $\frak a$-cofinite for all $i<n$, by the previous argument, $H_{\frak b}^i(M)$ is $\frak b$-cofinite for all $i<n$; and hence         
it follows from [NS, Theorem 3.7] that $\Hom_A(A/\frak b,H_{\frak b}^i(M))$ is finitely generated for all $i\leq n$. We observe that $\Hom_A(A/\frak a,\Gamma_{\frak b}(M))=\Hom_A(A/\frak a,\Gamma_{\frak a}(M))=\Hom_A(A/\frak a,M)$ is finitely generated. Furthermore, $\Hom_A(A/\frak a,H_{\frak b}^1(M))\cong\Hom_A(A/\frak a,\Hom_A(A/xA,H_{\frak b}^1(M))\cong \Hom_A(A/\frak b,H_{\frak b}^1(M))$ is finitely generated. Now, since $\Gamma_{\frak a}(M)$ is $\frak a$-cofinite, applying the functor $\Hom_A(A/\frak a,-)$ to the above exact sequence, we conclude that $\Hom_A(A/\frak a,H^1_{\frak a}(M))$ is finitely generated. By the argument mentioned in the beginning of the proof $\Hom_A(A/\frak a,H_{\frak a}^i(M))\cong\Hom_A(A/\frak a,H_{\frak b}^i(M))\cong\Hom_A(A/\frak b,H_{\frak b}^i(M))$ is finitely generated for all $2\leq i\geq n$.   
\end{proof}

It was proved in [NS, Theorems 3.3, 3.7] that if $\dim A/\frak a\leq 2$, then the conditions in \cref{three} are equivalent. In the following theorem we extend this result for $\dim A/\frak a\geq 3$, but by an additional assumption.
\medskip

\begin{Theorem}\label{tl}
Let $\dim A/\frak a=d\geq 3$ and let $\depth(\Ann(M),A/\frak a)\geq d-2$. If $\Ext_A^i(A/\frak a,M)$ is finitely generated for all $i\leq n+1$, then the conditions in \cref{three} are equivalent.
\end{Theorem}
\begin{proof}
 We proceed by induction on $d$. The case $d=3$ follows from \cref{three} and so we may assume that $d\geq 4$.
  Since $\depth(\Ann(M),A/\frak a)>0$, there exists an element $x\in\Ann_RM$ which is a non-zerodivisor on $A/\frak a$. Taking $\frak b=\frak a+xA$, we have $\dim A/\frak b=d-1$. Since $\Supp_AA/\frak b\subseteq \Supp_AA/\frak a$, it follows $\Ext_A^i(A/\frak b,M)$ is finitely generated for all $i\leq n+1$. If $H_{\frak a}^i(M)$ is $\frak a$-cofinite for each $i<n$, then $H_{\frak b}^i(M)$ is $\frak b$-cofinite for each $i<n$. The induction hypothesis implies that $\Hom_A(A/\frak b,H_{\frak b}^i(M))$ is finitely generated for all $i\leq n$. By the same reasoning in the proof of \cref{three}, we have $\Hom_A(A/\frak a,H_{\frak a}^i(M))$ is finitely generated for all $i\leq n$. Conversely, assume that $\Hom_A(A/\frak a,H_{\frak a}^i(M))$ is finitely generated for all $i\leq n$ and so $\Hom_A(A/\frak b,H_{\frak b}^i(M))$ is finitely generated for all $i\leq n$. The induction hypothesis implies that  $H_{\frak b}^i(M)$ is $\frak b$-cofinite for each $i<n$ and so $H_{\frak a}^i(M)$ is $\frak a$-cofinite for each $i<n$. 
\end{proof}

\medskip
\begin{Corollary}
Let $\frak b$ an ideal of $A$ and $\frak a=\Gamma_{\frak b}(A)$ such that $\dim A/\frak a=3$. Then $\Hom_A(A/\frak a,H_{\frak a}^i(A/\frak b))$ is finitely generated for all $i\leq n$ if and only if $H^i_{\frak a}(A/\frak b)$ is $\frak a$-cofinite for all $i<n$.
\end{Corollary}
\begin{proof}
As $\Gamma_{\frak b}(A/\frak a)=0$, we have $\depth(\frak b,A/\frak a)>0$. Now the assertion is obtained by using \cref{three}.  
\end{proof}

\medskip

\begin{Corollary}
Let $\frak b$ an ideal of $A$ and $\frak a=\Gamma_{\frak b}(A)$ such that $\dim A/\frak a=3$. Then $\Hom_A(A/\frak a,H_{\frak a}^i(A))$ is finitely generated for all $i\leq n$ if and only if $H^i_{\frak a}(A)$ is $\frak a$-cofinite for all $i<n$. In particular, $\Hom_A(A/\frak a,H_{\frak a}^1(A))$ is finitely generated.
\end{Corollary}
\begin{proof}
There exists a positive integer $t$ such that $\frak b^t\frak a=0$ and so $\Gamma_{\frak a}(\frak b^t)=\frak b^t$. Thus applying the functor $\Gamma_{\frak a}(-)$ to the exact sequence $0\To\frak b^t\To A\To A/\frak b^t\To 0$, we deduce that $H_{\frak a}^i(A)\cong H_{\frak a}^i(A/\frak b^t)$ for each $i> 0$. Since $\Gamma_{\frak b}(A/\frak a)=0$, we have $\depth(\frak b^t,A/\frak a)=\depth(\frak b,A/\frak a)>0$. If $\Hom_A(A/\frak a,H_{\frak a}^i(A))$ is finitely generated for all $i\leq n$, then $\Hom_A(A/\frak a,H_{\frak a}^i(A/\frak b^t))$ is finitely generated for all $i\leq n$. Now \cref{three} implies that $H_{\frak a}^i(A/\frak b^t)$ is $\frak a$-cofinite for all $i<n$; and hence $H_{\frak a}^i(A)$ is $\frak a$-cofinite for all $i<n$. Conversely, if $H_{\frak a}^i(A)$ is $\frak a$-cofinite for all $i<n$, then $H_{\frak a}^i(A/\frak b^t)$ is $\frak a$-cofinite for all $i<n$ and so using again \cref{three}, $\Hom_A(A/\frak a,H_{\frak a}^i(A/\frak b^t)$ is finitely generated for all $i\leq n$ so that $\Hom_A(A/\frak a,H_{\frak a}^i(A))$ is finitely generated for all $i\leq n$ 
\end{proof}

\medskip
\begin{Corollary}
Let $\frak p$ be a prime ideal of $A$ with $\dim A/\frak p=3$ and let $\frak b$ be an ideal of $A$ such that $\frak b\nsubseteq \frak p$. Then $\Hom_A(A/\frak p,H_{\frak p}^i(A/\frak b))$ is finitely generated for all $i\leq n$ if and only if  $H_{\frak p}^i(A/\frak b)$ is $\frak p$-cofinite for all $i<n$.
\end{Corollary}
\begin{proof}
Since $\frak b\nsubseteq \frak p$, we have $\depth(\frak b,A/\frak p)>0$ and so the the result follows from \cref{tl}.
\end{proof}

If $\dim A\geq 4$, then we have the following result.

\medskip

\begin{Proposition}
Let $\dim A=d\geq 4$ with $\depth(\Ann (M),A/(\frak a+\Gamma_{\frak a}(A)))\geq d-3$ and let $\Ext_A^i(A/\frak a,M)$ is finitely generated for all $i\leq n+1$. Then the conditions in \cref{three} are equivalent.
\end{Proposition}
\begin{proof}
We can choose an integer $t$ such that $(0:_A \frak a^t)=
\Gamma_{\frak a}(A)$. Put $\overline{A} = A/{\Gamma_{\frak a}(A)}$
and $\overline{M} = M/{(0:_M \frak a^t)}$ which is an
$\overline{A}$-module. Taking $\overline{\frak a}$ as the image of
$\frak a$ in $\overline{A}$, we have $\Gamma_{\overline{\frak
a}}({\overline{A}})=0$. Thus $\overline{\frak a}$ contains an
$\overline{A}$-regular element so that $\dim A/\frak a+\Gamma_{\frak
a}(A)={\rm dim}\overline{A}/\overline{\frak a}\leq d-1$. The
assumption on $M$ together with the fact that $ {\Supp} _{A}(A/\frak
a+\Gamma_{\frak a}(A))\subset {\Supp} _{A}(A/\frak a)$ and [DM,
Proposition 1] imply that ${\rm Ext}^i_A(A/\frak a+\Gamma_{\frak
a}(A),M)$ is finitely generated for all $i\leq n$. Since by the assumption $(0:_M\frak a)$ is finitely generated, it is clear that $(0:_M\frak a^t)$ is finitely generated and we have an exact sequence $$0\To (0:_M\frak a^t)\To \Gamma_{\frak a}(M)\To \Gamma_{\frak a}(\overline{M})\To 0$$ and the isomorphism $H_{\frak a}^i(M)\cong H_{\frak a}^i(\overline{M})$ for all $i>0$. In order to prove (i)$\Rightarrow$(ii), assume that $\Hom_A(A/\frak a,H_{\frak a}^i(M))$ is finitely generated for all $i\leq n$. Then in view of the previous argument and the independence theorem for local cohomology  $\Hom_A(A/\frak a+\Gamma_{\frak a}(A),H_{\frak a+\Gamma_{\frak a}(A)}^i(\overline{M}))$ is finitely generated for all $i\leq n$. It now follows from \cref{tl} that $H_{\frak a+\Gamma_{\frak a}(A)}(\overline{M})$ is $\frak a+\Gamma_{\frak a}(A)$-cofinite for all $i<n$; and hence using the change of ring principle [DM, Proposition 2], $H_{\frak a}^i(\overline{M})$ is $\frak a$-cofinite for all $i<n$. Consequently , the previous argument implies that $H_{\frak a}^i(M)$ is $\frak a$-cofinite for all $i<n$. (ii)$\Rightarrow$(i). Assume that $H_{\frak a}^i(M)$ is $\frak a$-cofinite for all $i<n$. By the same reasoning as mentioned before, we deduce that $H_{\frak a+\Gamma_{\frak a}(A)}^i(\overline{M})$ is $\frak a+\Gamma_{\frak a}(A)$-cofinite for all $i<n$. Now, using again \cref{tl}, we deduce that $\Hom_A(A/\frak a,H_{\frak a}^i(\overline{M}))\cong \Hom_A(A/\frak a+\Gamma_{\frak a}(A),H_{\frak a+\Gamma_{\frak a}(A)}^i(\overline{M}))$ is finitely generated for all $i\leq n$ and consequently the previous argument yields that $\Hom_A(A/\frak a,H_{\frak a}^i(M))$ is finitely generated for all $i\leq n$.
\end{proof}

\medskip







\begin{thebibliography}{Saz21}


\bibitem[BNS1]{BNS1} K. Bahmanpour, R. Naghipour, M. Sedghi, {\em Cofiniteness with respect to ideals of small dimension}, Algebr. Represent. Theory, \textbf{18} (2014), no.~2, 369--379.

 \bibitem[BNS2]{BNS2} K. Bahmanpour, R. Naghipour, M. Sedghi, {\em On  the category of cofinite modules which is abelian}, Proc. Amer. Math. Soc, {\bf 142} (2014), no. 4, 1101--1107.
 
 
\bibitem [BS]{BS} M. Brodmann, R.Y. Sharp, \emph{Local Cohomology: an Algebraic Introduction with Geometric Applications},
 Cambridge Univ. Press, Cambridge, UK (1998).
 
 \bibitem [DM]{DM} D. Delfino, T. Marley,
\emph{Cofinite modules of local cohomology}, J. Pure. Appl, Algebra,
\textbf{121(1)} (1997), 45--52.

\bibitem [H]{H} R. Hartshorne,
\emph{Affine duality and cofiniteness}, Invent. Math. \textbf{9}
(1970), 145--164.
\bibitem[K]{K} K-I. Kawasaki, {\em On a category of cofinite modules which is abelian}, Math. Z, {\bf 269} (2011), no. 1-2, 587--608.  


\bibitem[KS]{KS} M. Khazaei, R. Sazeedeh, {\em A criterion for cofiniteness of modules}, Rend. Sem. Mat. Padova, to appear.  

\bibitem [M1]{M1}\label{2005} L. Melkersson, \emph{Modules cofinite with respect to an ideal},
J. Algebra \textbf{285} (2005), 649--668.

\bibitem [M2]{M2}\label{2005} L. Melkersson, \emph{Cofiniteness with respect to ideals of dimesnion one},
J. Algebra \textbf{372} (2012), 459--462.

\bibitem [NS]{NS} M. Nazari, R. Sazeedeh,
\emph{Cofiniteness with respect to two ideals and local cohomology}, Algebr Represent Theory, {\bf 22} (2019), 375--385.

\bibitem [S]{S} R. Sazeedeh, \emph{Cofiniteness of Koszul cohomology modules}, arXiv:submit/4690830 [math.AC] 13 Jan 2023.
\end{thebibliography}
\end{document}